\documentclass[a4paper,12pt, reqno]{amsart}
\usepackage{hyperref}
\usepackage{a4wide}
\usepackage{amsthm}
\usepackage{amssymb}
\usepackage{palatino}
\usepackage{lipsum}
\usepackage{longtable}
\usepackage{enumitem}
\usepackage[utf8]{inputenc}
\def\ds{\displaystyle}

\usepackage{rotating}
\usepackage{tikz,lipsum,lmodern}
\usepackage[most]{tcolorbox}

\newtheorem{thm}{Theorem}[section]

\newtheorem{lem}[thm]{Lemma}

\theoremstyle{definition}

\newtheorem{rem}{Remark}

\theoremstyle{remark}

\numberwithin{equation}{section}

\newcommand{\legendre}[2]{\genfrac{(}{)}{}{}{#1}{#2}}

\begin{document}
\title[Some congruences for $(\ell, k)$ and $(\ell, k, r)$-regular partitions]{Some congruences for $(\ell, k)$ and $(\ell, k, r)$-regular partitions}
\author[T Kathiravan, K Srinivas and Usha K Sangale]{T Kathiravan$^{(1)}$, K Srinivas$^{(2)}$ and Usha K Sangale$^{(3)}$}
\address{$^{(1)}$Department of Mathematics, IIT Madras, Chennai, Tamilnadu, 600 036, India}
\address{$^{(2)}$Institute of Mathematical Sciences, HBNI,  CIT Campus, Taramani, Chennai, Tamilnadu 600 113, India}
 \address{$^{(3)}$SRTM University, Vishnupuri, Nanded, Maharashtra 431 606, India}
\email{$^{(1)}$kkathiravan98@gmail.com}
\email{$^{(2)}$srini@imsc.res.in}
\email{$^{(3)}$ushas073@gmail.com}


\begin{abstract} Let $b_{\ell, k}(n), b_{\ell, k, r}(n)$  count the number of $(\ell, k)$, $(\ell, k, r)$-regular partitions respectively.  In this paper we shall derive infinite families of congruences for $b_{\ell, k}(n)$ modulo $2$ when $ (\ell, k) = (3,8),  (4, 7)$, for $b_{\ell, k}(n)$ modulo $8$, modulo $9$ and modulo $12$ when $(\ell, k) = (4, 9)$ and $b_{\ell, k, r}(n)$ modulo $2$ when $(\ell, k, r) = (3, 5, 8)$. 
\end{abstract}

\maketitle

\section{Introduction}

We shall use the following standard notation throughout the paper: for complex numbers $a, q$ with $ \mid q \mid < 1 $, define 
\[ f_i := {(q^i; q^i)}_{\infty}, \ \  i = 1, 2, 3, \cdots,  \mathrm{where} \, \  {(a; q)}_{\infty} = \prod_{m=0}^{\infty} (1 - a q^m)
\]

Recall that a partition of a positive integer $n$ is a non-increasing sequence of positive integers whose sum is $n$. An $\ell$-regular partition is a partition in which none of the parts is divisible by $\ell$. Denote by $b_{\ell}(n)$ the number of $\ell$-regular partitions with the convention $b_{\ell}(0) = 1$. Then the generating function for $b_{\ell}(n)$ is given by
\begin{equation}\label{r-0}
 \sum_{n=0}^{\infty}b_{\ell}(n)q^n=\frac{f_{\ell}}{f_1}.
\end{equation}
Similarly, for coprime positive integers $\ell, k$ a partition of a positive integer $n$ is called $(\ell, k)$-regular partition  if none of the parts is divisible by $\ell$ and $k$. Let $b_{\ell, k}(n)$ denote the number of such partitions of $n$ with the notation that $b_{\ell, k}(0) = 1$. Then the generating function for $b_{\ell, k}(n)$ is given by
\begin{equation}\label{r}
 \sum_{n=0}^{\infty}b_{\ell, k}(n)q^n=\frac{f_{\ell} f_k}{f_1f_{\ell k}}.
 \end{equation}
 Generalizing the above concept to three variables, we say that for $ \ell, k, r$ positive integers with $ \textrm{gcd} ( \ell, k, r) =1$, a partition is an $ (\ell, k, r)$-regular partition if none of the parts is divisible by $\ell$, $k$ and $r$. Let $b_{\ell, k, r}(n)$ count the number of such partitions of $n$. Put $b_{\ell, k, r}(0) = 1$. Then the generating function for $b_{\ell, k, r} (n)$ is given by
 \begin{equation}\label{r2}
 \sum_{n=0}^{\infty}b_{\ell, k, r}(n)q^n=\frac{f_{\ell} f_k f_r}{f_1f_{\ell k} f_{\ell r}f_{kr} f_{\ell k r}}.
  \end{equation}

\medskip
There is an extensive literature on the congruences for $\ell$-regular partitions (see  for example \cite{Ahmed},  \cite{Ojah}, \cite{cui-gu}, \cite{Cui} and the references in these papers).
In \cite{Naika} several infinite families of congruences modulo $2$ for $b_{3,5}(n)$ have been proved. A similar study was carried out in \cite{Prasad} for $b_{3,5}(n)\pmod{2}$ where the partition was considered with distinct part.

\medskip

In this article, we shall derive some infinite families of congruences for $b_{\ell, k}(n)$ modulo $2$ when $ (\ell, k) = (3,8),  (4, 7)$, for $b_{\ell, k}(n) $ modulo $8$, modulo $9$ and modulo $12$ when $(\ell, k) = (4, 9)$ and $b_{\ell, k, r}(n)$ modulo $2$ when $(\ell, k, r) = (3, 5, 8)$. 

\medskip

More precisely, we prove the following congruences.

\medskip

  \begin{thm}\label{t3.8}
  Let $b(n)$ be defined by
  \begin{equation}\label{e-1-1}
\sum_{n=0}^{\infty}b(n)q^n=f_1^2f_3.
\end{equation}
Let $p \geq 5$ be a prime and let $\legendre{\star}{p}$ denote the Legendre symbol. Define
 \begin{equation}\label{e-1-2}
\omega (p) := b\big({{5(p^2-1)}/{24}}\big) + \legendre{-6}{p}  \legendre{-5(p^2-1)/{24}}{p}.
\end{equation}
\begin{enumerate}
\item[(i)]  If $\omega (p) \equiv0\pmod2$, then for $n, j \geq 0$, if $p \not | n$, we have
\begin{equation}\label{e-1-3}
b_{3,8}\big(2p^{4j+3}n+5p^{4j+4}/{12} +7/{12} \big) \equiv 0 \pmod2.
\end{equation}
\item[(ii)] If $\omega (p) \not\equiv0\pmod2$, then for $n, j \geq 0$,
\begin{enumerate}
\item[(a)] if $p \not | 24n + 5$, we have
\begin{equation}\label{e-1-31}
b_{3,8}\big(2p^{6j+2}n+5p^{6j+2}/{12} +7/{12} \big) \equiv 0 \pmod2,
\end{equation}
\item[(b)]
 if $p \not | n$, we have
\begin{equation}\label{e-1-4}
b_{3,8}\big(2p^{6j+5}n+5p^{6j+6}/{12} +7/{12} \big) \equiv 0 \pmod2.
\end{equation}
\end{enumerate}
\end{enumerate}
  \end{thm}
 
 \begin{rem}
 It must be mentioned that in \cite[Theorem 1.5]{xia} the authors considered \eqref{e-1-1} and derived infinite family of congruences modulo $3$ for Ramanujan's $\phi$-function.
 \end{rem}
 
 \begin{rem}
 By Mathematica, we see that $b(10) = 0$. Therefore, from \eqref{e-1-2}, putting $p=7$, we get 
 $ \omega(7) \equiv 1 \pmod{2} $ and from \eqref{j7-0}, we get $ \omega (7) \equiv  \left( \frac{n-3}{7} \right) \pmod{2}$. Both these conditions imply that $ n \not\equiv 3\pmod{7}$. Thus, taking $j=0$ in \eqref{e-1-31}, for example, we get 
 \[
 b_{3,8} (686 n + 98 r + 21 ) \equiv 0 \pmod{2},  \quad r\in \{0, 1, 2, 4, 5, 6\}.
 \]
 \end{rem}
  
  \medskip
 \begin{thm}\label{t4.7} For $n\geq0$ and $p \geq 3$ a prime, we have
 \begin{itemize}
 \item[(i)] $ b_{4, 7}\left(14(7n+j)+13\right) \equiv 0\pmod2, \ \ \mathrm{where} \ \  j\in\{1,2,3,4,5,6\}$.\\
 \item[(ii)] $b_{4, 7}\left(98p^2n+\frac{98p^2+6}{8}\right) \equiv b_{4, 7}(98n+13)\pmod2$.\\
\item[(iii)] $ b_{4, 7}\left(98p^2n+\frac{49p(p+6j)+3}{4}\right) \equiv 0\pmod2, \ \ \mathrm{where} \ \ j\in\{1, \cdots, p-1\}$.\\

\end{itemize}
 \end{thm}

\medskip

 \begin{thm}\label{t4.9} For $n\geq0$, we have
 \begin{enumerate}
\item[(i)]  $b_{4, 9}(8n+7) \equiv 0\pmod{12}$.\\
 \item[(ii)] $b_{4, 9}(16n+15) \equiv 0\pmod8$.\\
  \item[(iii)] $\sum_{n=0}^{\infty}b_{4, 9}(16n+7)q^n \equiv 4f^6_1f_3\pmod8$.\\
  \item[(iv)] $b_{4,9}(32n + 29) \equiv 5 b_{4,9} (16n + 15) \pmod{9}$.
  \end{enumerate}
  \end{thm}
  
  \medskip

  \begin{thm}\label{t4.9a} Let $a(n)$ be defined by
  \begin{equation}\label{h1-0}
\sum_{n=0}^{\infty}a(n)q^n=f^6_1f_3,
\end{equation}
  and let $p \geq 5$ be a prime. Define
  \begin{equation}\label{h5-1}
\omega (p) := a\big({{3(p^2-1)}/{8}}\big) + \legendre{-6}{p} p^2 \legendre{-3(p^2-1)/{8}}{p}.
\end{equation}
\begin{enumerate}
\item[(i)]  If $\omega (p) \equiv0\pmod8$, then for $n, k \geq 0$, if $p \not | n$, we have
\begin{equation}\label{h12-0}
b_{4,9}\big(16p^{4k+3}n+6p^{4k+4} +1 \big) \equiv 0 \pmod8.
\end{equation}
\item[(ii)] Let $ p\equiv1\pmod8$ and $\omega (p) \not\equiv0,2,4,6\pmod8$. Then for $n, k \geq 0$

\begin{enumerate}
\item[(a)] if $p \not| 8n +3$, we have 
\begin{equation}\label{h12-0-1}
b_{4,9}\big(16p^{6k+2}n+6p^{6k+2} +1 \big) \equiv 0 \pmod8, 
\end{equation}
\item[(b)] if $p \not | n$, we have
\begin{equation}\label{h12-1}
b_{4,9}\big(16p^{6k+5}n+6p^{6k+6} +1 \big) \equiv 0 \pmod8.
\end{equation}
\end{enumerate}
\end{enumerate}
  \end{thm}
   
 \begin{thm}\label{t3.5.8-1} Let $n\geq 0$. We have
 \begin{itemize}
 \item[(i)] $
  b_{3, 5, 8}(16(5n+j)+3) \equiv 0\pmod2, \quad j\in\{2, 4\},$\\
  \item[(ii)] 
 $ b_{3, 5, 8}(80(5n+j)+51) \equiv 0\pmod2, \quad j\in\{2, 4\}$,\\
 \item[(iii)]
 $ b_{3, 5, 8}(400(5n+j)+291) \equiv 0\pmod2, \quad j\in\{1, 3\}$,\\
 \item[(iv)] $
 b_{3, 5, 8}(2000n+1091) \equiv b_{3,5,8}(80n+51)\pmod2.\label {et3}$
 \end{itemize}
 \end{thm}

 \begin{thm}\label{t3.5.8-2} For $n \geq 0$, we have
 \begin{equation}
 b_{3,5,8}\left(16\cdot5^{2k+1}n+\frac{26\cdot5^{2k+1}+23}{3}\right)\equiv b_{3,5,8}(80n+51)\pmod2.
 \end{equation}
 \end{thm}

 \section{Preliminaries}
 We recall that for $\mid ab\mid<1$, Ramanujan's general theta function $f(a,b)$ is defined as
\begin{equation} \label{al}
f(a,b)=\ds\sum^\infty_{n=-\infty}a^{n(n+1)/2}b^{n(n-1)/2}.
\end{equation}

\medskip

Now from Jacobi's triple product identity \cite[Entry 19, p. 35]{berndt2012ramanujan} 
\begin{equation}
f(a,b)=(-a;ab)_\infty(-b;ab)_\infty(ab;ab)_\infty,
\end{equation}
 it follows that (see \cite[Entry 22, p.36]{berndt2012ramanujan})
\begin{equation}
\varphi(q):=f(q,q)=1+2\ds\sum^\infty_{n=1}q^{n^2}=(-q;q^2)^2_\infty(q^2;q^2)_\infty=\frac{f^5_2}{f^2_1f^2_4},
\end{equation}
\begin{equation}
\psi(q):=f(q;q^3)=\ds\sum^\infty_{n=0}q^{n(n+1)/2}=\frac{(q^2;q^2)_\infty}{(q;q^2)_\infty}=\frac{f^2_2}{f_1},
\end{equation}
and
\begin{equation} \label{bl}
f(-q):=f(-q,-q^2)=\ds\sum^\infty_{n=-\infty}(-1)^nq^{n(3n-1)/2}=(q;q)_\infty=f_1.
\end{equation}
Now we shall state few lemmas which will be required for the proof of our results.
 
\medskip

 \begin{lem}\label{l2} \normalfont{\cite[Theorem 2.2]{cui-gu}} Let $p\geq5$ be a prime such that
\begin{equation}\nonumber
\frac{\pm p-1}{6}:=\begin{cases}
      \ds\frac{p-1}{6}, & \text{if}\ p\equiv1\pmod6, \\
      \ds\frac{-p-1}{6}, & \text{if}\ p\equiv-1\pmod6.
    \end{cases}
\end{equation}
With the notations as defined in \eqref{al} and \eqref{bl}, we have
\[
f_1 =\ds\sum^\frac{p-1}{2}_{\substack{k=-\frac{p-1}{2}\\k\neq\frac{\pm p-1}{6}}}(-1)^kq^\frac{3k^2+k}{2}
f\left(-q^\frac{3p^2+(6k+1)p}{2},-q^\frac{3p^2-(6k+1)p}{2}\right)
+(-1)^\frac{\pm p-1}{6}q^\frac{p^2-1}{24} f_{p^2}.
\]
Furthermore, if $-\frac{p-1}{2}\leq k\leq\frac{p-1}{2},  k\neq\frac{\pm p-1}{6}$, then
$\frac{3k^2+k}{2}\not\equiv\frac{p^2-1}{24}\pmod p$.
\end{lem}

 \medskip
\begin{lem}\label{l1}\normalfont {\cite[Lemma 2.3]{Ahmed}} If $p\geq3$ is prime, then
\[
f_1^3 = \sum^{p-1}_{\substack{{k=0}\\{k\neq\frac{p-1}{2}}}}(-1)^kq^\frac{k(k+1)}{2}
\sum^\infty_{n=0}(-1)^n(2pn+2k+1)q^{pn\cdot\frac{pn+2k+1}{2}}
+p(-1)^\frac{p-1}{2}q^\frac{p^2-1}{8} f_{p^2}^3.
\]
Furthermore, if $k\neq\frac{p-1}{2}, 0\leq k\leq p-1,$ then
$\frac{k^2+k}{2}\not\equiv\frac{p^2-1}{8}\pmod p.$
\end{lem}

\medskip

 \begin{lem} {\textnormal{\cite[Theorem 4.6]{Ojah}}} Define the generalized partition function $p_{c^sd^t}(n)$ by
 \begin{equation}\label{gen-p-1}
 \ds\sum^\infty_{n=0}p_{c^sd^t}(n)q^n=\frac{1}{f_c^sf_{d}^t}.
\end{equation}
Then 
\begin{eqnarray}
\sum^\infty_{n=0}p_{3^15^1}(2n+1)q^n&=&q\frac{f^2_2f^2_{30}}{f^2_3f^2_5f_1f_{15}} \label{b8}, \\  
\sum^\infty_{n=0}p_{1^1{15}^1}(2n+1)q^n&=&q\frac{f^2_6f^2_{10}}{f^2_1f_3f_5f^2_{15}} \label{b4}.
\end{eqnarray}
  \end{lem}

\medskip

\begin{lem}{\normalfont{\cite[Theorem 4.1]{Ojah}}} Define the generalized partition function $p_{a^x b^y c^z} (n) $ by
\begin{equation}\label{ab0}
\sum_{n=0}^{\infty} p_{a^x b^y c^z d^u}(n) q^n := \frac{1}{f_a^x f_b^y f_c^z f_d^u}. 
\end{equation}
Then
\begin{equation}\label{b2}
\ds\sum^\infty_{n=0}p_{1^13^15^115^1}(2n+1)q^n=\frac{f_2f_6f_{10}f_{30}}{f^2_1f^2_3f^2_5f^2_{15}}+2q\frac{f^2_2f^2_6f^2_{10}f^2_{30}}{f^3_1f^3_3f^3_5f^3_{15}}.
\end{equation}
\end{lem}

\medskip

\begin{lem} The following $2$-dissections hold:
 \begin{equation}\label{l7}
 \ds\frac{1}{f^2_1}=\frac{f^5_8}{f^5_2f^2_{16}}+2q\frac{f^2_4f^2_{16}}{f_2^5f_8},
 \end{equation}
 \begin{equation}\label{l6}
   \frac{1}{f^4_1} = \frac{f^{14}_4}{f^{14}_2f^4_8}+4q\frac{f^2_4f^4_8}{f^{10}_2}.
 \end{equation}
\end{lem}
See \cite[Entry 25, p. 40]{berndt2012ramanujan}
\medskip

\begin{lem} The following $2$-dissections hold:
 \begin{equation}\label{A1}
 \frac{f_3}{f_1}=\frac{f_4f_6f_{16}f^2_{24}}{f^2_2f_8f_{12}f_{48}}+q\frac{f_6f^2_8f_{48}}{f^2_2f_{16}f_{24}},
 \end{equation}
 \begin{equation}\label{l3}
 \frac{f_5}{f_1}=\frac{f_8f^2_{20}}{f^2_2f_{40}}+q\frac{f^3_4f_{10}f_{40}}{f^3_2f_8f_{20}},
 \end{equation}
 \begin{equation}\label{l4}
\frac{f_9}{f_1}=\frac{f^3_{12}f_{18}}{f^2_2f_6f_{36}}+q\frac{f^2_4f_6f_{36}}{f^3_2f_{12}}.
 \end{equation}
 \end{lem}
 The proofs of \eqref{A1}, \eqref{l3} and \eqref{l4} can be seen in \cite{XY}, \cite{hirschhorn2010elementary} and \cite{Xia}, respectively.
 \medskip
 
 \begin{lem} The following $2$-dissection holds:
 \begin{equation}\label{l5}
  \ds\frac{f_3}{f_1^3}=\frac{f_4^6f_6^3}{f_2^9 f_{12}^2}+3q\frac{f_4^2f_6 f_{12}^2}{f_2^7}.
 \end{equation}
 \end{lem}
 The proof can be found in  Baruah and Ojah \cite{Ojah} and Xia and Yao \cite{xia-yao}.
 
 \medskip
 
 \begin{lem} The following 2-dissection holds:
 \begin{equation}\label{l2}
 \ds\frac{1}{f_1f_7}=\frac{f_{16}^2f_{56}^5}{f^2_2f_8f_{14}^2f_{28}^2f_{112}^2}+q\frac{f^2_4f^2_{28}}{f^3_2f^3_{14}}+q^6\frac{f^5_8f^2_{112}}{f^2_2f^2_4f^2_{14}f^2_{16}f_{56}}.
 \end{equation}
 \end{lem}
 See  \cite[p.315]{berndt2012ramanujan}.
 
\medskip
The following result of M. Newman \cite[Theorem 3]{Newman} will play a crucial role in the proof of our theorems, therefore we shall quote is as a lemma. 
Following the notations of Newman's paper, we shall let $p,q$ denote distinct primes, let $r, s$ be integers such that $r, s \neq 0, r \not\equiv s\pmod{2}.$ Set
\begin{equation}\label{new-1}
\phi (\tau) = \prod {( 1 - x^n)}^r  {( 1 - x^{nq})}^s  = \sum c(n) x^n,
\end{equation}
 $ \quad \quad \varepsilon =  \frac{1}{2} (r +s),  t = ( r + sq )/{24}, \Delta = t (p^2 - 1)/{24}, \theta = { (-1)}^{\frac{1}{2} - \varepsilon} 2 q^s, \legendre{\star}{p}$ is the well-known Legendre symbol. Then the result is as follows.
\begin{lem}\label{newman}
With the notations defined as above, the coefficients $c(n)$ of $\phi(\tau)$ satisfy
\begin{equation}\label{new-2}
c(np^2 + \Delta ) - \gamma_n c(n) + p^{2\varepsilon - 2} c\big((n - \Delta)/{p^2}\big) = 0,
\end{equation}
where
\begin{equation}\label{new-3}
\gamma_n = p^{2\varepsilon - 2} c - \legendre{\theta}{p} p^{\varepsilon - 3/2} \legendre{n - \Delta}{p}.
\end{equation}
\end{lem}

\medskip

\noindent
We shall also require the following result which is known as Ramanujan's $5$-dissection of Euler's product. Ramanujan stated this formula without a proof. See \cite[p.85]{Hirsch} for a proof.

\medskip

\begin{lem} We have
\begin{equation}\label{h}
f_1 = f_{25} \big(R^{-1} - q - q^2 R \big), \ \textrm{where} \ R = \frac{ {(q^5 ; q^{25})}_{\infty}{(q^{20} ; q^{25})}_{\infty}}{{ {(q^{10} ; q^{25})}_{\infty}}{(q^{15} ; q^{25})}_{\infty}}.
\end{equation}

\end{lem}

\section{Proof of theorems}

\subsection{Proof of theorem \ref{t3.8}} By definition
\begin{equation} \label{j1}
\sum_{n=0}^{\infty}b_{3, 8}(n)q^n=\frac{f_3f_8}{f_1f_{24}}.
\end{equation}
Substituting \eqref{A1} into \eqref{j1} and extracting terms involving $q^{2n+1}$ from both sides,  we get
\begin{equation} \label{j1-a}
 \sum_{n=0}^{\infty}b_{3, 8}(2n+1)q^{2n+1} =q  \frac{f_6 f_8^3 f_{48}}{f_2^2 f_{16} f_{24}^2}.
 \end{equation}
Cancelling $q$ from both sides, then replacing $q^2$ by $q$, and using the binomial relation
\begin{equation}\label{bin-1}
f_k^2 \equiv f_{2k} \pmod{2},
\end{equation}
we get
  \begin{equation} \label{j3}
  \sum_{n=0}^{\infty}b_{3, 8}(2n+1)q^n \equiv \frac{f_3f^3_4f_{24}}{f^2_1f_8f^2_{12}}\equiv f_1^2 f_3  \pmod{2}.
 \end{equation}
We now write 
\begin{equation} \label{1j3}
\sum_{n=0}^{\infty}b(n)q^n= f^2_1f_3.
\end{equation}
Putting $r=2, q=3$ and $s = 1$ in \eqref{new-1}, we have by Lemma \ref{newman}, for any $n \geq 0$
\begin{equation}\label{j4}
 b \big(p^2n+{5(p^2-1)}/{24}\big)=\gamma_n b(n)-p \ b\big( (n- 5(p^2-1)/{24}) / {p^2}\big),
 \end{equation}
 where
\begin{equation}\label{j5}
\gamma_n=p c - \legendre{-6}{p} \legendre{ n- 5(p^2-1)/{24}}{p}
\end{equation} 
and $c$ is a constant.

\medskip

\noindent
Setting $n=0$ in \eqref{j4} gives us
$$ b \big({5(p^2-1)}/{24}\big) = \gamma_0,
$$
since $b(n)=0,$ if $n$ is not an integer. On the other hand, from \eqref{j5}, we have
$$
\gamma_0=p c - \legendre{-6}{p} \legendre{ - 5(p^2-1)/{24}}{p}.
$$
This gives
$$ 
p c =  b \big({5(p^2-1)}/{24}\big) + \legendre{-6}{p} \legendre{ - 5(p^2-1)/{24}}{p}.
$$ 
Now the equation \eqref{j4} can be rewritten as 
\begin{equation}\label{j7}
 b \big(p^2n+{5(p^2-1)}/{24}\big) = \big( \omega(p) - \legendre{-6}{p} \legendre{ n- 5(p^2-1)/{24}}{p}  \big) b(n) - p b \big( (n- 5(p^2-1)/{24}) / {p^2}\big),
\end{equation}
where
\begin{equation*}
\omega(p) = b \big( {5(p^2-1)}/{24}\big) + \legendre{-6}{p} \legendre{- 5(p^2-1)/{24}}{p}.
\end{equation*}
Note that, if 
\begin{equation}\label{j7-0}
 \omega(p) = \legendre{-6}{p} \legendre{n - 5(p^2-1)/{24}}{p},
 \end{equation}
  then by \eqref{j7}, we get
\begin{equation}\label{j7-1}
 b \big(p^2n+{5(p^2-1)}/{24}\big) \equiv p b \big( (n- 5(p^2-1)/{24}) / {p^2}\big) \pmod{2}
\end{equation}
Now, if $\omega(p)\not \equiv 0\pmod{2}$, then we note that $p \not| {24n + 5}$. Therefore, for such primes we see that 
$ (n- 5(p^2-1)/{24}) / {p^2}$ is not an integer. Thus, we get
\begin{equation}\label{j7-2}
 b \big(p^2n+{5(p^2-1)}/{24}\big) \equiv 0 \pmod{2}.
 \end{equation}
\noindent
Now, replacing $n$ by $pn+ {5(p^2-1)}/{24}$ in \eqref{j7}, we get
\begin{equation}\label{j11}
 b \big(p^3n+{5(p^4-1)}/{24}\big) = \omega(p)  b \big(pn+{5(p^2-1)}/{24}\big) - p b (n/p).
\end{equation}
If $ \omega (p) \equiv 0\pmod{2}$, then \eqref{j11} yields
\begin{equation}\label{j12}
 b \big(p^3n+{5(p^4-1)}/{24}\big) \equiv  p b (n/p) \pmod{2}.
\end{equation}
Replacing $n$ by $pn$ in \eqref{j12} gives
\begin{equation}\label{j13}
 b \big(p^4n+{5(p^4-1)}/{24}\big) \equiv  p b (n) \pmod{2}.
\end{equation}
By induction, this gives
\begin{equation}\label{j14}
 b \big(p^{4j} n+{5(p^{4j}-1)}/{24}\big) \equiv  p^{4j} b (n) \pmod{2}.
\end{equation}
If $ p \not| n$, then \eqref{j12} yields
\begin{equation}\label{j15}
 b \big(p^3n+{5(p^4-1)}/{24}\big) \equiv  0 \pmod{2}.
\end{equation}
Thus replacing $n$ by $\big(p^3n+{5(p^4-1)}/{24}\big)$ in \eqref{j14} and using \eqref{j15}, we get
\begin{equation}\label{j16}
 b \big(p^{4j+3} n+{5(p^{4j+4}-1)}/{24}\big) \equiv  0 \pmod{2}.
\end{equation}
This settles the proof of part (i) of the theorem by observing that $ b(n) \equiv b_{3,8} (2n+1) \pmod{2}$.

\medskip

In order to prove part (ii), we replace $n$ by $p^2n+\big( {5 p (p^2-1)}/{24}\big)$ in \eqref{j11}
\begin{eqnarray}\label{j17}
b \big(p^5n+{5(p^6-1)}/{24}\big) & \equiv  &\omega (p) b \big(p^3n+{5(p^4-1)}/{24}\big) - p \  b\big(pn+{5(p^2-1)}/{24}\big).\nonumber \\
& \equiv & \big[ \omega^2 (p) - p \big] b \big(pn+{5(p^2-1)}/{24}\big) - p \omega (p) \ b(n/p).
\end{eqnarray}
Now, if $\omega (p) \not\equiv 0 \pmod{2}$, then $\omega^2 (p) - p \equiv 0 \pmod2$, and therefore \eqref{j17} becomes
\begin{equation}\label{j18}
b \big(p^5n+{5(p^6-1)}/{24}\big) \equiv b(n/p) \pmod{2}.
\end{equation}
Replacing $n$ by $pn$, \eqref{j18} yields
\begin{equation}\label{j19}
b \big(p^6n+{5(p^6-1)}/{24}\big) \equiv  b(n) \pmod{2}.
\end{equation}
 By induction, we can write the above equation as
\begin{equation}\label{j20}
b \big(p^{6j}n+{5(p^{6j}-1)}/{24}\big) \equiv  b(n) \pmod{2}.
\end{equation}
Replacing $n$ by $ p^{2}n+{5(p^{2}-1)}/{24}$ in \eqref{j20} and using \eqref{j7-2}, we obtain
\begin{equation}\label{j20-1}
b \big(p^{6j+2}n+{5(p^{6j+2}-1)}/{24}\big) \equiv  0 \pmod{2}.
\end{equation}

\noindent
Observe that if $p \not| n$, then $b(n/p) =0$. Thus \eqref{j18} yields
\begin{equation}\label{j21}
b \big(p^5n+{5(p^6-1)}/{24}\big) \equiv 0 \pmod{2}.
\end{equation}
Replacing $n$ by $p^5n+{5(p^6-1)}/{24}$ in \eqref{j20} and using \eqref{j21}, we obtain
\begin{equation}\label{j22}
b \big(p^{6j+5}n+{5(p^{6j+6}-1)}/{24}\big) \equiv 0 \pmod{2}.
\end{equation}
The proof of part (ii) now follows from  \eqref{j20-1} and \eqref{j22} by observing that $b(n) \equiv b_{3,8}(2n+1) \pmod{2}$.

\subsection{Proof of theorem \ref{t4.7}}

We take $\ell = 4, k = 7$ in \eqref{r} to get
\begin{equation}\label{d1}
\sum_{n=0}^{\infty} b_{4,7} (n) q^n = \frac{f_4 f_7}{f_1 f_{28}} \equiv \frac{f_1^3 f_7}{f_{28}}\pmod{2}.
\end{equation}
Here, we have used $f_4 \equiv f_1^4\pmod{2}$ which follows by binomial theorem. Now
replacing the expression $f_1^3$ in \eqref{d1} by Lemma \ref{l1} with $p=7$ and then extracting those terms involving $q^{7n+6}$, we obtain
\begin{equation}\label{d1-1}
 \sum_{n=0}^{\infty}b_{4, 7}(7n+6)q^{7n+6} \equiv q^6 \frac{f_7 f_{49}^3}{f_{28}}  \pmod{2}.
 \end{equation}
Cancelling $q^6$ from both sides, then replacing $q^7$ by $q$ in \eqref{d1-1}, we get
\begin{equation}\label{d2}
 \sum_{n=0}^{\infty}b_{4, 7}(7n+6)q^n\equiv \frac{f_1f_7^3 }{f_{4}} \equiv \frac{f_1f_7 f_{14} }{f_{4}}\pmod2.
 \end{equation}
Observe that by \eqref{l2}
\begin{equation}\label{d2-1}
f_1 f_7 \equiv \frac{f_2 f_{14}}{f_1 f_7} \equiv f_2 f_{14} \Big\{ \frac{f_{16}^2f_{56}^5}{f^2_2f_8f_{14}^2f_{28}^2f_{112}^2}+q\frac{f^2_4f^2_{28}}{f^3_2f^3_{14}}+q^6\frac{f^5_8f^2_{112}}{f^2_2f^2_4f^2_{14}f^2_{16}f_{56}}\Big\} \pmod{2}.
\end{equation}
Now,  substituting \eqref{d2-1} in \eqref{d2}, and extracting only those terms involving $q^{2n +1}$ on both sides, we get
 \begin{equation}\label{d2-2}
 \sum_{n=0}^{\infty}b_{4, 7}(14n+13)q^{2n+1}\equiv q \frac{f_{4} f_{28}^2}{f_2^2 f_{14}}   \pmod2.
 \end{equation}
Cancelling $q$ from both sides, using the binomial relation $f_{ 2^k 7} \equiv f_7^{2^k } \pmod{2}$ and then 
replacing $q^2$ by $q$ in \eqref{d2-2}, we obtain
  \begin{equation}\label{d3}
 \sum_{n=0}^{\infty}b_{4, 7}(14n+13)q^n\equiv f^3_7\pmod2.
 \end{equation}
Note that on the right hand side of \eqref{d3}, the powers of $q$ are multiples of $7$,  thus part (i) follows.
 
 \medskip
 
 Now to establish the remaining parts of the theorem, we extract only those terms in \eqref{d3} in which the power of $q$ is congruent to $0$ modulo $7$, replace $q^7$ by $q$ and obtain
  \begin{equation}\label{d4}
 \sum_{n=0}^{\infty}b_{4, 7}(98n+13)q^n\equiv f^3_1\pmod2.
 \end{equation}
 
 \medskip
 
 Now, substituting Lemma \ref{l1} in \eqref{d4}, extracting those terms in which the power of $q$ is congruent to $(p^2-1)/8$ modulo $p$, we obtain
 \begin{equation}\label{d5}
  \sum_{n=0}^{\infty}b_{4, 7}\big(98 (pn+\frac{p^2-1}{8} )+ 13\big)q^{pn+ \frac{p^2-1}{8}} \equiv  q^{\frac{p^2-1}{8}}f_{p^2}^3 \pmod2.
 \end{equation}
 Cancelling $q^{\frac{p^2-1}{8}}$, then replacing $q^p$ by $q$, \eqref{d5} gives 
 \begin{equation}\label{d6}
  \sum_{n=0}^{\infty}b_{4, 7}\big(98 (pn+\frac{p^2-1}{8} )+ 13\big)q^n \equiv  f_{p}^3 \pmod2.
 \end{equation}
Note that the right hand side of \eqref{d6} contains only terms in $q$ whose powers are congruent to $0$ modulo $p$. Thus ,
 \begin{equation}\label{d7}
  \sum_{n=0}^{\infty}b_{4, 7}\big(98 (p^2n+\frac{p^2-1}{8} )+ 13\big)q^{pn} \equiv  f_{p}^3 \pmod2.
 \end{equation}
and 
 \begin{equation}\label{d8}
  \sum_{n=0}^{\infty}b_{4, 7}\big(98 \big( p (pn +j) +\frac{p^2-1}{8} \big)+ 13\big)q^n \equiv  0  \pmod2, \quad j= 1, 2, \cdots, p-1.
 \end{equation}
Proof of part (ii) of the theorem follows on taking $q^p$ as $q$ in \eqref{d7}, then comparing with \eqref{d4} and part (iii) follows from \eqref{d8}.


\subsection{Proof of theorem \ref{t4.9}}

Setting $\ell=4$ and $k=9$ in \eqref{r}, we have
 \begin{equation}\label{f1}
 \sum_{n=0}^{\infty}b_{4, 9}(n)q^n=\frac{f_4f_9}{f_1f_{36}}.
 \end{equation}
Now, substituting \eqref{l4} in \eqref{f1} and extracting those terms in which the power of $q$ is congruent to $1$ modulo $2$, we get
 \begin{equation}\label{f1-1}
 \sum_{n=0}^{\infty}b_{4, 9}(2n+1)q^{2n+1}= q \frac{f_4^3f_6}{f_2^3f_{12}}.
 \end{equation}
Replacing $q^2$ by $q$ in the above equation, we get
 \begin{equation}\label{f2}
 \sum_{n=0}^{\infty}b_{4, 9}(2n+1)q^n=\frac{f^3_2f_3}{f^3_1f_6}.
 \end{equation}
Substituting \eqref{l5} into \eqref{f2}, we have  
\begin{equation}\label{f3}
 \sum_{n=0}^{\infty}b_{4, 9}(2n+1)q^n=\frac{f^3_2}{f_6}\left( \frac{f_4^6f_6^3}{f_2^9 f_{12}^2}+3q\frac{f_4^2f_6 f_{12}^2}{f_2^7}\right).
\end{equation}
Now, extracting those terms in which the power of $q$ is congruent to $1$ modulo $2$, replacing $q^2$ by $q$, we have  
\begin{equation}\label{f4}
\sum_{n=0}^{\infty}b_{4, 9}(4n+3)q^n=3\frac{f^2_2f^2_6}{f^4_1}.
\end{equation}
Again, substituting \eqref{l6} into \eqref{f4}, and extracting those terms in which the power of $q$ is congruent to $1$ modulo $2$, replacing $q^2$ by $q$, we have
 \begin{equation}\label{f5}
\sum_{n=0}^{\infty}b_{4, 9}(8n+7)q^n=12\frac{f^2_2f^2_3f^4_4}{f^8_1}.
\end{equation}
Thus the proof of (i) follows from \eqref{f5}. 

\medskip
\noindent
Now we use the congruence $f_1^{2^k} \equiv f_2^{2^{k-1}} \pmod{2^k}$ ( see eqn. (3.7), \cite{yao}) to obtain
$f_1^8 \equiv f_2^4\pmod{8}$ and use the congruence $f_4^4 \equiv f_2^8\pmod{8}$ in \eqref{f5} to obtain
\begin{equation}\label{f5a}
\sum_{n=0}^{\infty}b_{4, 9}(8n+7)q^n\equiv4f_2^6f_6\pmod8.
\end{equation}
From \eqref{f5a} extracting the powers of $q$ which are congruent to $1$ modulo $2$, we obtain the proof of (ii) and extracting the powers of $q$ which are congruent to $0$ modulo $2$, replacing $q^2$ by $q$, the proof of part (iii) is established.

\medskip

To prove part (v) which involves congruences modulo $9$, we note that 
$f_3^2 \equiv f_1^6\pmod{3}$, thus \eqref{f5} can be written as
\begin{equation}\label{f7}
\sum_{n=0}^{\infty}b_{4, 9}(8n+7)q^n\equiv3\frac{f^2_2f^4_4}{f^2_1}\pmod9.
\end{equation}
Replacing  $1/{f_1^2}$ in \eqref{f7} with \eqref{l7} and extracting terms which involve odd powers of $q$, we obtain
\begin{equation}\label{f7a}
\sum_{n=0}^{\infty}b_{4, 9}(16n+15)q^{2n+1} \equiv 6 q \frac{f_4^6 f_{16}^2}{f_2^3 f_8}\pmod9.
\end{equation}
Cancelling $q$ from both sides, replacing $q^2$ by $q$, \eqref{f7a} gives
\begin{eqnarray}\label{f7b}
\sum_{n=0}^{\infty}b_{4, 9}(16n+15)q^{n} &\equiv & 6  \frac{f_2^6 f_{8}^2}{f_1^3 f_4}\equiv 6  \frac{f_6^2 f_{8}^2}{f_3 f_4}\pmod9, \nonumber\\
& \equiv & 6 \psi (q^3) \psi (q^4)\pmod9,
\end{eqnarray}
where $\psi(q) := {f_2^2}/{f_1}$. 

\medskip

Now, we first extract terms involving $q^{2n}$ from \eqref{f2}, which gives us
\begin{equation}\label{g1}
 \sum_{n=0}^{\infty}b_{4, 9}(4n+1)q^{2n} =\frac{f_4^6f_6^2}{f_2^6f_{12}^2}.
 \end{equation}
Now replace $q^2$ by $q$ in \eqref{g1} to get
\begin{equation}\label{g2}
 \sum_{n=0}^{\infty}b_{4, 9}(4n+1)q^{n} =\frac{f_2^6f_3^2}{f_1^6f_{6}^2}.
 \end{equation}
Substituting \eqref{l5} for $\frac{f_3^2}{f_1^6}$ in \eqref{g2}, extracting terms containing $q^{2n+1} $on both sides, we obtain
\begin{equation}\label{g3}
 \sum_{n=0}^{\infty}b_{4, 9}(8n+5)q^{2n+1} =6 q \frac{f_4^8 f_6^2}{f_2^{10}}.
 \end{equation}
Cancelling $q$ on both sides, replacing $q^2$ by $q$, this yields
\begin{equation}\label{g4}
 \sum_{n=0}^{\infty}b_{4, 9}(8n+5)q^{n} =6 \frac{f_2^8 f_3^2}{f_1^{10}}.
 \end{equation}
Using the congruence $f_3^2 \equiv f_1^6 \pmod{3}$, we obtain from \eqref{g4}
\begin{equation}\label{g5}
 \sum_{n=0}^{\infty}b_{4, 9}(8n+5)q^{n} \equiv 6 \frac{f_2^8}{f_1^{4}}\pmod{9}.
 \end{equation}
Substituting \eqref{l6} in \eqref{g5}, extracting terms containing $q^{2n+1}$ from both sides, we obtain
\begin{equation}\label{g6}
 \sum_{n=0}^{\infty}b_{4, 9}(16 n+13)q^{2n+1} \equiv 6q \frac{f_4^2f_8^4}{f_2^2}\pmod{9}.
 \end{equation}
Cancelling $q$ from both sides of \eqref{g6}, replacing $q^2$ by $q$, we get
\begin{equation}\label{g7}
 \sum_{n=0}^{\infty}b_{4, 9}(16 n+13)q^{n} \equiv 6 \frac{f_2^2f_4^4}{f_1^2}\pmod{9}.
 \end{equation}
Now, in place of  $1/f_1^2$ in \eqref{g7}, we substitute the expression \eqref{l7}, then extract terms on both sides containing $q^{2n+1}$, we get
\begin{equation}\label{g8}
 \sum_{n=0}^{\infty}b_{4, 9}(32 n+29)q^{2n+1} \equiv 3 q\frac{f_4^6f_{16}^2}{f_2^3 f_8}\pmod{9}.
 \end{equation}
Cancelling $q$ from both sides, replacing $q^2$ by $q$, \eqref{g8} gives
\begin{eqnarray}\label{g9}
 \sum_{n=0}^{\infty}b_{4, 9}(32 n+29)q^{n} &\equiv& 3 \frac{f_2^6f_{8}^2}{f_1^3 f_4} \equiv 3 \frac{f_6^2f_{8}^2}{f_3 f_4}\pmod{9},\nonumber\\
&\equiv& 3 \psi(q^3) \psi(q^4).
 \end{eqnarray}
Thus from \eqref{f7b} and \eqref{g9}, it follows that
$$
b_{4, 9}(32 n+29) \equiv 5  b_{4, 9}(16 n+15)\pmod{9}.
$$
\qed

\subsection{Proof of theorem \ref{t4.9a}} 
\medskip

Let us set
\begin{equation}\label{h1}
\sum_{n=0}^{\infty}a(n)q^n=f^6_1f_3.
\end{equation}
It follows from the statement (iii) of theorem \ref{t4.9} that
\begin{equation}\label{part-3}
\sum_{n=0}^{\infty}b_{4, 9}(16n+7)q^n \equiv 4\sum_{n=0}^{\infty}a(n)q^n \pmod8.
\end{equation}

\medskip

\noindent
Taking $r = 6, s = 1$ and $q = 3$ in Lemma\eqref{newman}, we see that for any $n \geq 0$, we have
\begin{equation}\label{h2}
 a\big(p^2n+{3(p^2-1)}/{8}\big)=\gamma_n a(n)-p^5 a\big( \big({n-{3(p^2-1)}/{8}}\big)/{p^2}\big),
 \end{equation}
 where
\begin{equation}\label{h3}
\ds\gamma_n = p^5c-\legendre{-6}{p} \ p^2 \  \legendre{{n - 3(p^2-1)/8}}{p}
\end{equation} 
and $c$ is a constant.

\medskip

\noindent
Note that $a\big( {{3(p^2-1)}/{8{p^2}}}\big)=0$, since $ {{3(p^2-1)}/{8{p^2}}}$ is not an integer. Also $a(0) = 1$. Thus, 
from \eqref{h2}, we get $ \gamma_0 = a\big({{3(p^2-1)}/{8}}\big)$.

\medskip

\noindent
Again from \eqref{h3}, it follows that
\begin{equation}\nonumber
\gamma_0 = p^5c-\legendre{-6}{p} p^2 \legendre{-3(p^2-1)/{8}}{p}.
\end{equation}
Therefore,
\begin{equation}\label{h4}
p^5c = a\big({{3(p^2-1)}/{8}}\big) + \legendre{-6}{p} \ p^2 \ \legendre{-3(p^2-1)/{8}}{p}.
\end{equation}
Thus the equation \eqref{h2} can be rewritten as
\begin{multline}\label{h5}
a\big(p^2n+{3(p^2-1)}/{8}\big) \\ = \left[\omega (p) - p^2 \  \legendre{{-6n + 9(p^2-1)/4}}{p} \right] a(n) - p^5 a\big( \big({n-{3(p^2-1)}/{8}}\big)/{p^2}\big),
\end{multline}
where
\begin{equation*}
\omega (p) = a\big({{3(p^2-1)}/{8}}\big) + \legendre{-6}{p} p^2 \legendre{-3(p^2-1)/{8}}{p}.
\end{equation*}

\medskip
Note that, if 
\begin{equation}\label{h5-0}
 \omega(p) = p^2 \left(\frac{-6n+9(p^2-1)/4}{p}  \right)
 \end{equation}
  then by \eqref{h5}, we get
\begin{equation}\label{h5-1}
 a \big(p^2n+{3(p^2-1)}/{8}\big) \equiv  -p^5 a \big( (n- 3(p^2-1)/{8}) / {p^2}\big) \pmod{8}
\end{equation}
Now, if $\omega(p)\not \equiv 0\pmod{8}$, then we note that $p \not| {8n + 3}$. Therefore, for such primes we see that 
$ (n- 3(p^2-1)/{8}) / {p^2}$ is not an integer. Thus, we get
\begin{equation}\label{h5-2}
 a \big(p^2n+{3(p^2-1)}/{8}\big) \equiv 0 \pmod{8}.
 \end{equation}

\noindent
Replacing $n$ by $pn+{3(p^2-1)}/{8}$ in \eqref{h5}, we obtain
\begin{equation}\label{h9}
a\big(p^3n+{3(p^4-1)}/{8}\big) = \omega (p) a(pn+{3(p^2-1)}/{8}) - p^5 a(n/p).
\end{equation}
Note that if $\omega (p) \equiv0\pmod8$, we get
\begin{equation}\label{h10}
a\big(p^3n+{3(p^4-1)}/{8}\big) \equiv - p^5 a(n/p)\pmod{8}.
\end{equation}
Replacing $n$ by $pn$, we get
\begin{equation}\label{h11}
a\big(p^4n+{3(p^4-1)}/{8}\big) \equiv -p^5 a(n)\pmod8.
\end{equation}
By mathematical induction, it gives
\begin{equation}\label{h12}
a\big(p^{4j}n+{3(p^{4j}-1)}/{8}\big) \equiv - p^{5j} a(n)\pmod8.
\end{equation}

\medskip

\noindent
If $p \not| n$, then \eqref{h10} gives
\begin{equation}\label{h13}
a\big(p^3n+{3(p^4-1)}/{8}\big) \equiv 0\pmod{8}.
\end{equation}
Replacing $n$ by $p^3n+{3(p^4-1)}/{8}$ in \eqref{h12} and using \eqref{h13}, we get
\begin{equation}\label{h14}
a\big(p^{4j+3}n+3(p^{4j+4}-1)/{8} \big) \equiv 0\pmod{8}.
\end{equation}

\noindent
Again replacing $n$ by $p^2n+{3p(p^2-1)}/{8}$ in \eqref{h9}, we get
\begin{equation}\label{h15}
a\big(p^5n+{3(p^6-1)}/{8}\big) = \omega (p) a(p^3n+{3(p^4-1)}/{8}) - p^5 a(pn+{3(p^2-1)}/{8}).
\end{equation}
Combining \eqref{h9} and \eqref{h15}, we obtain
\begin{equation}\label{h16}
a\big(p^5n+{3(p^6-1)}/{8}\big) = (\omega^2 (p) - p^5) a\big(pn+{3(p^2-1)}/{8}\big) - \omega p^5 a(n/p).
\end{equation}
If $\omega(p)\not\equiv0,2,4,6\pmod8$ and $p\equiv1\pmod8$, then
$\omega^2 (p) - p^5\equiv0\pmod8$,
therefore \eqref{h16} becomes
\begin{equation}\label{h17}
a\big(p^5n+{3(p^6-1)}/{8}\big)\equiv - \omega(p)  p^5 a(n/p)\pmod{8}.
\end{equation}
Now we replace $n$ by $pn$ in \eqref{h17} to obtain
\begin{equation}\label{h18}
a\big(p^6n+{3(p^6-1)}/{8}\big)\equiv - \omega(p) p^5 a(n)\pmod{8}.
\end{equation}
Using mathematical induction, \eqref{h18} gives
\begin{equation}\label{h19}
a\big(p^{6j}n+{3(p^{6j}-1)}/{8}\big)\equiv - \omega (p) p^{5j} a(n)\pmod{8}.
\end{equation}

Replacing $n$ by $ p^{2}n+{3(p^{2}-1)}/{8}$ in \eqref{h19} and using \eqref{h5-2}, we obtain
\begin{equation}\label{h19-1}
 a \big(p^{6j+2}n+{3(p^{6j+2}-1)}/{8}\big) \equiv  0 \pmod{8}.
\end{equation}

Observe that if $p \not | n$, then \eqref{h17} gives
\begin{equation}\label{h20}
a\big(p^5n+{3(p^6-1)}/{8}\big)\equiv 0\pmod{8}.
\end{equation}
Therefore, replacing $n$ by $p^5n+{3(p^6-1)}/{8}$ in \eqref{h19} and then using \eqref{h20}, we get
\begin{equation}\label{h21}
a\big(p^{6j+5}n+{3(p^{6j+6}-1)}/{8}\big)\equiv 0\pmod{8}.
\end{equation}
As $b_{4, 9}(16n+7) \equiv a(n) \pmod{8}$, the theorem follows from \eqref{h14}, \eqref{h19-1} and \eqref{h21} respectively.

\medskip

\subsection{Proof of theorem \ref{t3.5.8-1}} Setting $\ell=3, k=5$ and $r=8$ in \eqref{r2}, we obtain
 \begin{equation}\label{e1}
 \sum_{n=0}^{\infty}b_{3, 5, 8}(n)q^n=\frac{f_3f_5f_8}{f_1f_{15}f_{24}f_{40}f_{120}}.
 \end{equation}
 Using the relation
 \[
 f_{24} \equiv f_6^4 \equiv f_6^3 f_6 \equiv f_6^3 f_3^2\pmod{2}, \quad
 f_{40} \equiv f_{10}^4 \equiv f_{10} f_{10}^3 \equiv f_5^2 f_{10}^3\pmod{2},
 \]
we can write \eqref{e1}  as 
 \begin{equation}\label{e2}
 \sum_{n=0}^{\infty}b_{3, 5, 8}(n)q^n\equiv\frac{f_8}{f^3_{6}f^3_{10}f_{120}f_1f_3f_5f_{15}}\pmod2.
 \end{equation}
As $\quad \frac{1}{f_1f_3f_5f_{15}} := \sum^\infty_{n=0}p_{1^13^15^115^1}(n)q^n$, \eqref{e2} can be rewritten as
\begin{equation}\label{e2a}
 \sum_{n=0}^{\infty}b_{3, 5, 8}(n) q^n \equiv \frac{f_8}{f^3_{6}f^3_{10}f_{120}} \sum^\infty_{n=0}p_{1^13^15^115^1}(n)q^n\pmod{2}.
\end{equation}
In \eqref{e2a}, extracting the powers of $q$ which are congruent to $1$ modulo $2$, then cancelling $q$ from both sides and replacing $q^2$ by $q$, we obtain
\begin{equation}\label{e2b}
 \sum_{n=0}^{\infty}b_{3, 5, 8}(2n+1)q^n\equiv\frac{f_4}{f_{3}^3f_{5}^3f_{60}} \sum^\infty_{n=0}p_{1^13^15^115^1}(2n+1)q^n \pmod2.
 \end{equation}
 It follows from \eqref{b2} that $ \sum^\infty_{n=0}p_{1^13^15^115^1}(2n+1)q^n \equiv 1\pmod{2}$, therefore, \eqref{e2b} yields
 \begin{equation}\label{e3}
 \sum_{n=0}^{\infty}b_{3, 5, 8}(2n+1)q^n\equiv  \frac{f_4}{f_{3}^3f_{5}^3f_{60}} \equiv \frac{f_4}{f_3f_5f_{6}f_{10}f_{60}}\pmod2.
 \end{equation}
 Note that $\frac{1}{f_3f_5} :=  \sum^\infty_{n=0}p_{3^1 5^1}(n)q^n.$ Thus, we can write \eqref{e3} as
  \begin{equation}\label{e3a}
 \sum_{n=0}^{\infty}b_{3, 5, 8}(2n+1)q^n \equiv  \frac{f_4}{f_{6}f_{10}f_{60}}  \sum^\infty_{n=0}p_{3^1 5^1}(n)q^n\pmod2.
 \end{equation}
 
 Extracting from both sides of \eqref{e3a} only those terms involving $q^{2n+1}$, cancelling $q$ from both sides, then replacing $q^2$ by $q$, we get
  \begin{equation}\label{e3b}
 \sum_{n=0}^{\infty}b_{3, 5, 8}(4n+3)q^n  \equiv  \frac{f_2}{f_{3}f_{5}f_{30}}  \sum^\infty_{n=0}p_{3^1 5^1}(2n+1)q^n\pmod2.
 \end{equation}
 Replacing the sum on the right hand side of \eqref{e3b} with \eqref{b8}, we obtain
 \begin{equation}\label{e3c}
 \sum_{n=0}^{\infty}b_{3, 5, 8}(4n+3)q^n \equiv  q \frac{f_2^3 f_{30}}{f_1f_3f_5f_6f_{10}f_{15}} \pmod{2}.
 \end{equation}
 Again we use  $\quad \frac{1}{f_1f_3f_5f_{15}} := \sum^\infty_{n=0}p_{1^13^15^115^1}(n)q^n$ in \eqref{e3c} to obtain
 \begin{align}\label{e3d}
 \sum_{n=0}^{\infty}b_{3, 5, 8}(4n+3)q^n & \equiv  q \frac{f_2^3 f_{30}}{f_6f_{10}} \sum^\infty_{n=0}p_{1^13^15^115^1}(n)q^n \pmod{2}, \nonumber\\
 & \equiv q \frac{f_2^3 f_{30}}{f_6f_{10}} \big\{ \sum^\infty_{n=0} p_{1^13^15^115^1}(2n)q^{2n} + \sum^\infty_{n=0}p_{1^13^15^115^1}(2n+1)q^{2n+1}\big\}\pmod{2}.
 \end{align}
 From \eqref{e3d} we first extract only those terms on both sides in which the power of $q$ is congruent to $0$ modulo $2$, then replace $q^2$ by $q$  and then use the relation $ \sum^\infty_{n=0}p_{1^13^15^115^1}(2n+1)q^n \equiv 1\pmod{2}$ to obtain
 \begin{align}
 \sum_{n=0}^{\infty}& b_{3, 5, 8}(8n+3)q^n  \equiv q\frac{f_1^3f_{15}}{f_3f_5}  \equiv q\frac{f_1^4f_{15}^2}{f_1f_3f_5f_{15}} \pmod{2}, \nonumber \\
 & \equiv q\frac{f_4f_{30}}{f_1f_3f_5f_{15}} \equiv q f_4f_{30} \sum^\infty_{n=0}p_{1^13^15^115^1}(n)q^n \pmod2. \label{e5}
 \end{align}
 Now, comparing even powers of $q$, replacing $q^2$ by $q$ and then by \eqref{b2}, \eqref{e5} can be written as 
  \begin{align}
 \sum_{n=0}^{\infty} b_{3, 5, 8}(16n+3)q^n  \equiv q f_2f_{15} \sum^\infty_{n=0}p_{1^13^15^115^1}(2n+1)q^n \equiv q f_2f_{15} \pmod2.
 \label{e6}
 \end{align}
 We replace $q$ by $q^2$ in \eqref{h} to obtain
 \begin{equation}\label{h1}
f_2 = f_{50} \big(S^{-1} - q^2 - q^4 S \big), \ \textrm{where} \ S = \frac{ {(q^{10} ; q^{50})}_{\infty}{(q^{40} ; q^{50})}_{\infty}}{{ {(q^{20} ; q^{50})}_{\infty}}{(q^{30} ; q^{50})}_{\infty}}.
\end{equation}
 Now substituting \eqref{h1} in \eqref{e6}, we get
  \begin{equation}\label{e8}
 \sum_{n=0}^{\infty}b_{3, 5, 8}(16n+3)q^n\equiv f_{15}f_{50}\left(qS^{-1}-q^3-q^5S\right)\pmod2.
 \end{equation}
 Observe that the right hand side of \eqref{e8} doesn't contain terms of the form $q^{5n+2}$ or $q^{5n+4}$. Thus
  \begin{equation}\label{e8a}
 \sum_{n=0}^{\infty}b_{3, 5, 8}\big(16(5n+j) + 3\big)q^n\equiv 0 \pmod2, \ \textrm{for} \ j = 2, 4.
 \end{equation}
Extracting terms of the form $q^{5n+3}$ from both sides of \eqref{e8}, then replacing $q^5$ by $q$, we get
\begin{equation}\label{e9}
 \sum_{n=0}^{\infty}b_{3, 5, 8}(80n+51)q^n\equiv f_3f_{10}\pmod2.
 \end{equation}
Replacing $q$ by $q^3$ in \eqref{h} and substituting it in \eqref{e9}, we get
 \begin{equation}\label{e10}
 \sum_{n=0}^{\infty}b_{3, 5, 8}(80n+51)q^n\equiv f_{10}f_{75}\left(S^{-1}_1-q^3-q^6S_1\right)\pmod2, 
 \end{equation}
 where $ S_1 = \frac{ {(q^{15} ; q^{75})}_{\infty}{(q^{60} ; q^{75})}_{\infty}}{{ {(q^{30} ; q^{75})}_{\infty}}{(q^{45} ; q^{75})}_{\infty}}$.
 
 \noindent
 We note that on the right hand side of \eqref{e10}, the terms containing $q^{5n +2}$ or $q^{5n +4}$ are absent, thus we infer
 \begin{equation}\label{e10a}
 \sum_{n=0}^{\infty}b_{3, 5, 8}\big(80 (5n +j) + 51\big)q^n\equiv 0 \pmod2,  \ \textrm{for} \ j = 2, 4.
 \end{equation}
 From \eqref{e10}, extracting only those terms in which the power of $q$ is congruent to $3$ modulo $5$ and then replacing $q^5$ by $q$, we obtain
 \begin{equation}\label{e11}
 \sum_{n=0}^{\infty}b_{3, 5, 8}(400n+291)q^n\equiv f_{2}f_{15}\pmod2.
 \end{equation}
 Substituting \eqref{h1} in \eqref{e10}, we have
  \begin{equation}\label{e12}
 \sum_{n=0}^{\infty}b_{3, 5, 8}(400n+291)q^n\equiv f_{15}f_{50}\left(S^{-1}-q^2-q^4S\right)\pmod2.
 \end{equation}
Note that the terms containing $q^{5n+1}$ or $q^{5n+3}$ are absent on the right hand side of \eqref{e12}.
 Thus, we get
 \begin{equation}\label{e12a}
 \sum_{n=0}^{\infty}b_{3, 5, 8}\big(400(5n+j) + 291\big)q^n\equiv 0 \pmod2,  \ \textrm{for} \ j = 1, 3.
 \end{equation}
 Comparing the terms containing $q^{5n+2}$ on both sides of \eqref{e12}, then replacing $q^5$ by $q$, we get
  \begin{equation}\label{e13}
 \sum_{n=0}^{\infty}b_{3, 5, 8}(2000n+1091)q^n\equiv f_{3}f_{10} \equiv \sum_{n=0}^{\infty}b_{3, 5, 8}(80n+51)q^n\ \pmod2, 
 \end{equation}
 by \eqref{e9}. 
 
 \medskip
 
 \noindent
 The proof of Theorem \ref{t3.5.8-1} follows from \eqref{e8a}, \eqref{e10a}, \eqref{e12a} and \eqref{e13}.
 
\medskip
 
\subsection{Proof of theorem \ref{t3.5.8-2}}

For simplicity, we write \eqref{e9} as 
\begin{equation}\label{e14}
 \sum_{n=0}^{\infty}B(n)q^n\equiv f_3f_{10}\pmod2.
 \end{equation}
Then, \eqref{e11} becomes
\begin{equation}\label{e15}
 \sum_{n=0}^{\infty}B(5n + 3)q^n\equiv f_{2}f_{15}\pmod2.
 \end{equation}
 Extracting terms containing $q^{5n+2}$ from both sides of \eqref{e15}, repeating the arguments that led to \eqref{e13}, we get
\begin{equation*}
 \sum_{n=0}^{\infty}B \big(5 (5n+2) +3\big)q^n\equiv f_3f_{10}\pmod2
 \end{equation*}
or 
\begin{equation}\label{e16}
 \sum_{n=0}^{\infty}B(5^2n+ 2\times 5 +3)q^n\equiv f_3f_{10}\pmod2
 \end{equation}
 Extracting terms containing $q^{5n+3}$ from both sides of \eqref{e16}, repeating the arguments that led to \eqref{e11}, we get
\begin{equation*}
 \sum_{n=0}^{\infty}B \big(5 (5n+3) + 2 \times 5 +3\big)q^n\equiv f_2f_{15}\pmod2
 \end{equation*}
or 
\begin{equation}\label{e17}
 \sum_{n=0}^{\infty}B(5^3n+ 3\times 5^2 + 2\times 5 +3)q^n\equiv f_2f_{15}\pmod2
 \end{equation}
Repeating the above arguments one more time, we get
\begin{equation*}
 \sum_{n=0}^{\infty}B( 5^4n + 2 \times 5^3 + 3 \times 5^2 + 2\times 5 +3)q^n\equiv f_3f_{10}\pmod2
 \end{equation*}
Arguing inductively, at $m^{\textrm{th}}$ step, we get
\begin{equation*}
 \sum_{n=0}^{\infty}B( 5^{2m} n + 2 \times 5^{2m-1} + 3 \times 5^{2m-2} + \cdots  + 2\times 5 +3)q^n\equiv f_3f_{10}\pmod2
 \end{equation*}
Simplifying the above expression, we obtain
\begin{equation}\label{e18}
 \sum_{n=0}^{\infty}B \big( 5^{2m} n + 13 (5^{2m}-1)/24 \big) q^n\equiv f_3f_{10}\pmod2
 \end{equation}
 As $B(n) := b_{3, 5, 8}(80n+51)$, the theorem now follows from \eqref{e18}.

\medskip

\noindent
\textbf{Acknowledgements.} The first author acknowledges with thanks the financial support received from SERB NPDF grant no. PDF/2020/001755.


\hfill


\begin{thebibliography}{99}
 
 \bibitem{Ahmed}
 Zakir Ahmed and N. D. Baruah, New congruences for $\ell$-regular partitions for $\ell\in\{5, 6, 7, 49\}$, Ramanujan J., $40$ $(2016)$, $649-668$.
 
 \bibitem{Ojah}
N. D. Baruah and K. K. Ojah, Analogues of Ramanujan's partition identities and congruences arising from his theta functions and modular equations, Ramanujan J., $28$ $(2012)$, $385 - 407$.

\bibitem{berndt2012ramanujan}
B. C. Berndt, Ramanujan’s Notebooks, Part III, Springer, New York, 1991.

\bibitem{cui-gu}
S. P. Cui and N. S. S. Gu, Arithmetic properties of $\ell$-regular partitions, Adv. Appl. Math., $51$ $(2013)$, $507 - 523$.


 \bibitem{Cui}
S. P. Cui and N. S. S. Gu, Congruences for $\ell$-regular partitions and bipartitions, Rocky Mountain J. Math., $50$ $(2020)$, $513-526$.


\bibitem{Hirsch}
M. D. Hirschhorn, The powers of $q$, Developments in Mathematics, $49$, Springer International Publishing AG 2017.

\bibitem{hirschhorn2010elementary}
M. D. Hirschhorn and J. A. Sellers, Elementary proofs of parity results for $5$-regular partitions, Bull. Aust. Math. Soc., $81$ $(2010)$,  $58–63$.

\bibitem{Naika} M. S. Mahadeva Naika, B. Hemanthkumar and H. S. Sumanth Bharadwaj, Congruences modulo $2$ for certain partition functions, Bull. Aust. Math. Soc., $93$ $(2016)$, $400-409$.

\bibitem{Newman}
M. Newman, Modular forms whose coefficients possess multiplicative properties, $II$,  Ann. Math., $75$ $(1962)$, $242–250$. 

\bibitem{Prasad}
M. Prasad and K. V. Prasad, On $(\ell, m)$-regular partitions with distinct parts, Ramanujan J., $46$ $ (2018)$, $19 - 27$.



\bibitem{yao}
O. X. M. Yao, Congruences modulo $16, 32$, and $64$ for Andrew's singular overpartitions, Ramanujan J., $43$ $(2017)$, $215-225$.

\bibitem{xia} E. X. W. Xia, New congruence properties for Ramanujan's $\phi$ function, Proc. AMS., $149$ $(2021)$, $4985-4999$.

\bibitem{Xia}
E. X. W. Xia and O .X. M. Yao, Some modular relations for the Göllnitz–Gordon functions by an even–odd method, J. Math. Anal. Appl., $387$ $(2012)$, $126-138$.

\bibitem{xia-yao}
E. X. W. Xia and O .X. M. Yao, Analogues of Ramanujan's partition identities, Ramanujan J., $31$ $(2013)$, $373-396$.

\bibitem{XY}
E. X. W. Xia and O .X. M. Yao, New Ramanujan-like congruences modulo powers of 2 and 3 for overpartitions, J. Number Theory, $133$ $(2013)$, $1932-1949$.
\end{thebibliography}
\end{document}